\documentclass[11pt,draft]{amsart}
\usepackage{amssymb,amsmath}
\newtheorem{theo}{Theorem}[section]
\newtheorem{lema}{Lemma}[section]
\newtheorem{prop}{Proposition}[section]

\newcommand{\bq}{\begin{equation}}
\newcommand{\eq}{\end{equation}}
\newcommand{\ba}{\begin{array}}
\newcommand{\ea}{\end{array}}

\newcommand{\refe}[1]{(\ref{#1})}
\newcommand{\half}{\frac{1}{2}}

\newcommand{\dst}{\displaystyle}

\newcommand\be{\begin{enumerate}}
\newcommand\ee{\end{enumerate}}
\newcommand\bi{\begin{itemize}}
\newcommand\ei{\end{itemize}}

\newcommand{\pe}[2]{\langle#1,#2\rangle}

\newcommand\XX[1]{\mathbb{#1}}
\newcommand{\dq}{\mbox{$\Delta^{(1)}$}}

\newcommand\bd{\begin{definicion}{\bf }}
\newcommand\ed{\end{definicion}}
\newcommand\bl{\begin{lema}{\bf }}
\newcommand\el{\end{lema}}
\newcommand\bp{\begin{prop}{\bf }}
\newcommand\ep{\end{prop}}
\newcommand\bt{\begin{theo}{\bf }}
\newcommand\et{\end{theo}}
\newcommand\bdm{\begin{proof}}
\newcommand\edm{\end{proof}}
\newcommand\bn{\begin{nota}{\bf }}
\newcommand\en{\end{nota}}
\newcommand\bc{\begin{corolary}{\bf }}
\newcommand\ec{\end{corolary}}

\newtheorem{nota}{Remark}[section]

\newtheorem{definicion}{Definition}[section]
\newtheorem{corolary}{Corollary}[section]
\renewenvironment{proof}{\noindent{\bf Proof:\ }}
{\hfill \fbox \par\vspace{2ex}}

\newfont{\got}{eufm10 scaled \magstep1}
\newcommand{\g}[1]{\mbox{\got #1}}

\begin{document}

\title{Second structure relation for $q$-semiclassical
polynomials of the Hahn Tableau}
\author{R.S. Costas-Santos, F. Marcell\'{a}n}
\address{Roberto S. Costas-Santos. Dpto. de Matem\'{a}ticas,
Universidad Carlos III de Madrid, Spain.}
\email{rcostas@math.uc3m.es}

\address{Francisco Marcell\'{a}n. Dpto. de Matem\'{a}ticas,
Universidad Carlos III de Madrid, Spain.}

\email{pacomarc@ing.uc3m.es}
\thanks{}

\thanks{}
\date{\today}

\dedicatory{}

\commby{}
\begin{abstract}
$q$-Classical orthogonal polynomials of the $q$-Hahn
tableau are characterized from their orthogonality condition
and by a first and a second structure relation. Unfortunately,
for the $q$-semiclassical orthogonal polynomials (a
generalization of the classical ones) we find only
in the literature the first structure relation. In this paper,
a second structure relation is deduced. In particular, by means of a general
finite-type relation between a $q$-semiclassical polynomial
sequence and the sequence of its $q$-differences such a
structure relation is obtained.
\end{abstract}
\maketitle
\section{Introduction} \noindent
The $q$-Classical orthogonal polynomial sequences (Big
$q$-Jacobi, $q$-Laguerre, Al-Salam Carlitz I, $q$-Charlier,
etc.) are characterized by the property that the sequence of
its monic $q$-difference polynomials is, again, orthogonal
(Hahn's property, see \cite{hahn}). In fact, the
$q$-difference operator is a particular case of the Hahn
operator which is defined as follows
$$
L_{q,\omega}(f)(x)=\frac{f(qx+\omega)-f(x)}{(q-1)x+\omega},\qquad
\omega \in \XX C, \, q\in \XX C,\ |q|\ne 1.
$$
In the sequel, we are going to work with $q$-semiclassical
orthogonal polynomials and $q$-classical polynomials of the
Hahn Tableau, hence we will consider the $q$-linear lattice
$x(s)$, i.e. $x(s+1)=qx(s)+\omega$.
Therefore, for the sake of convenience we will denote
$\dq\equiv L_{q,\omega}$. Notice that for $q=1$ we get the
forward difference operator $\Delta$. In such a case,
when $w\to 0$ we recover the standard semiclassical orthogonal
polynomials \cite{Mar1}.
\\
Taking into account the role  of such families of
$q$-polynomials in the analysis of hypergeometric
$q$-difference equations resulting from physical
problems as the $q$-Schr\"{o}dinger equation,
$q$-harmonic oscillators, the connection and the linearization
problems among others there is an increasing interest to study them.
Moreover, the connection between the
representation theory of quantum algebras and the $q$-orthogonal
polynomials is well known (see \cite{renarv} and references therein).
\\
We also find many different approaches to the subject in the literature.
For instance, the functional equation (the so-called Pearson equation)
satisfied by the corresponding moment functionals allows an efficient
study of some properties of $q$-classical polynomials \cite{renmed},
\cite{kher}, \cite{makh}, \cite{medrenmar}.
However, the $q$-classical sequences of orthogonal polynomials
$\{\text{C}_n\}_{n\geq0}$ can also be characterized taking into
account its orthogonality as well as one of the two following
difference equations, the so-called structure relations.
\be
\item[$\bullet$]{\it First structure relation} \cite{als}, \cite{koedie},
\cite{nsu}
\bq \label{fistre}
\Phi(s)C_n^{[1]}(s)=\;
\sum_{\nu=n}^{n+t}\lambda_{n,\nu} \text{C}_{\nu}(s),\;
n\geq 0,\; \quad \lambda_{n,n}\ne 0,\; n\geq 0,
\eq
where $\Phi$ is a polynomial with $\deg \Phi=t\leq 2$ and
$C_n^{[1]}(s):=[n+1]^{-1}\dq\text{C}_{n+1}(s)$,
being
$$
[n]:=(q^n-1)/(q-1),\quad n\ge 0.
$$
\item[$\bullet\bullet$]{\it Second structure relation}
\cite{med,medrenmar}
\bq \label{sestre}
\text{C}_n(s)=\sum_{\nu=n-t}^{n}\theta_{n,\nu}
C_\nu^{[1]}(s),\quad \ n\geq t, \ 0\le t\le 2, \quad
\theta_{n,n}=1, \ n\ge t.
\eq
\ee
The $q$-classical orthogonal polynomials were introduced by W. Hahn
\cite{hahn} and also analyzed in \cite{als}. The generalization of this families leads to
$q$-semiclassical orthogonal polynomials which were introduced
by P. Maroni and extensively studied in the last decade by
himself, L. Kheriji, J. C. Medem, and others (see \cite{kher,med}).
\\
For $q$-classical orthogonal polynomial sequences, which are
$q$-semiclassical of class zero, the structure
relations  \refe{fistre} and \refe{sestre} become
$$
\ba{l}
\phi(s)L_{q,\omega} P_n(s)=\widetilde{\alpha}_n
P_{n+1}(s)+\widetilde {\beta}_n P_{n}(s)+\widetilde{\gamma}_n
P_{n-1}(s),\qquad \widetilde \gamma_n\ne 0,\hfill \\[0.4cm]
\sigma(s)L_{1/q,\omega/q} P_n(s)=\widehat{\alpha}_n
P_{n+1}(s) +\widehat{\beta}_n P_{n}(s)+\widehat{\gamma}_n
P_{n-1}(s), \qquad \widehat \gamma_n\ne 0, \\[0.4cm] P_n(s)=P_n^{[1]}(s)+\delta_n P_{n-1}^{[1]}
(s) +\epsilon_n P_{n-2}^{[1]}(s).\hfill
\ea
$$
In particular, in Table \ref{tab1} we describe  these parameters
for some families of $q$-classical orthogonal polynomials.
\begin{table}[!hbt] \label{tab1}
$$
\ba{l}
\hline \\[-0.3cm]
(A_1) \quad \mbox{Big $q$-Jacobi} \quad \widehat P_n(x;a,b,c;q)
\qquad x\equiv x(s)=q^{s}\\[0.2cm]
P^{[1]}_{n}(x;a,b,c;q)=q^{-n}\widehat P_n(qx;aq,bq,cq;q)
\\[0.2cm]
\phi(x)=aq(x-1)(bx-c) \qquad \sigma(x)=q^{-1}(x-aq)(x-cq)
\\[0.2cm]
\widehat{\alpha}_n= abq[n]\qquad \widetilde\alpha_n=q^{-n}[n]
\\[0.2cm]
\dst \widehat{\beta}_n=-aq[n](1-abq^{n+1})\frac{c+ab^2q^{2n+1}+
b(1-cq^n-cq^{n+1} -aq^n(1+q-cq^{n+1}))}{(1-abq^{2n})
(1-abq^{2n+2})} \\[0.5cm]
\dst \widetilde{\beta}_n=q[n](1-abq^{n+1})
\frac{c+a^2bq^{2n+1}+ a(1-cq^n-cq^{n+1}-bq^n(1+q-cq^{n+1}))}
{(1-abq^{2n}) (1-abq^{2n+2})}\\[0.5cm]
\dst \widehat{\gamma}_n=aq[n]\frac{(1-aq^n)(1-bq^n)(1-abq^n)
(c-abq^n)(1-cq^n)(1-abq^{n+1})}{(1-abq^{2n})^2(1-abq^{2n-1})
(1-abq^{2n+1})}\\[0.5cm]
\dst \widetilde \gamma_n=q^n \widehat \gamma_n \qquad
\delta_n=-\frac{q^n(1-q)}{1-abq^{n+1}}\widehat \beta_n \qquad
\epsilon_n=abq^{2n}\frac{(1-q^{n-1})(1-q)}{(1-abq^{n})
(1-abq^{n+1})}\widehat \gamma_n\\[0.4cm]
\hline \\[-0.3cm]
(A_2) \quad \mbox{$q$-Laguerre} \quad \widehat L_n^{(\alpha)}
(x;q) \qquad x\equiv x(s)=q^{s}\\[0.2cm]
L_{n}^{[1]\,(\alpha)}(x;q)=q^{-n}\widehat L_n^{(\alpha+1)}
(qx;q)\\[0.2cm]
\phi(x)=ax(x+1) \qquad \sigma(x)=q^{-1}x\\[0.2cm]
\dst \widehat{\alpha}_n= a[n]\qquad
\widehat{\beta}_n=q^{-2n-1}[n](1+q-aq^{n+1})\qquad
\widehat{\gamma}_n=a^{-1}q^{1-4n}[n](1-aq^n)\\[0.3cm]
\dst \widetilde\alpha_n=0\qquad
\widetilde{\beta}_n=q^{-n}[n] \qquad \widetilde \gamma_n=
a^{-1}q^{1-3n}(1-aq^n) \\[0.3cm]
\delta_n=a^{-1}(1-q)\widehat \beta_n \qquad
\epsilon_n=a^{-1}(1-q^{n-1})(1-q)\widehat\gamma_n\\[0.3cm]
\hline \\[-0.3cm]
(A_3) \quad \mbox{Al-Salam Carlitz I} \quad \widehat
U_n^{(a)}(x;q) \qquad x\equiv x(s)=q^{s} \\[0.1cm]
U_{n}^{[1]\,(a)}(x;q)=\widehat U_n^{(a)}(x;q)\\[0.1cm]
\phi(x)=a \quad \sigma(x)=(1-x)(a-x) \quad
\widetilde\alpha_n=q^{1-n}[n]\quad \widetilde{\beta}_n=
q(1+a)[n] \quad \widetilde \gamma_n= aq^n[n] \\[0.2cm]
\hline\\[-0.3cm]
(A_4) \quad \mbox{$q$-Charlier} \quad \widehat
C_n(q^{-s};a;q)\\[0.2cm]
C_n^{[1]}(q^{-s};a;q)=\widehat C_n(q^{-s};aq^{-1};q)\\[0.2cm]
\phi(x)=x(x-1) \qquad \sigma(x)=q^{-1}ax \\[0.2cm]
\dst \widehat{\alpha}_n= [n]\qquad
\widehat{\beta}_n=q^{-2n-1}[n](a+aq+q^{n+1})\qquad
\widehat{\gamma}_n=aq^{1-4n}[n](a+q^n)\\[0.3cm]
\dst \widetilde\alpha_n=0\qquad
\widetilde{\beta}_n=aq^{-n}[n] \qquad \widetilde \gamma_n=
q^n \widehat{\gamma}_n \qquad \delta_n=(1-q)\widehat
\beta_n \qquad \epsilon_n=(1-q^{n-1})(1-q)\widehat\gamma_n
\\[0.3cm]  \hline
\ea
$$
\caption{Some families of $q$-polynomials of the Hahn Tableau}
\end{table}
\\
The first structure relation for the $q$-semiclassical
orthogonal polynomials was established
(see \cite{kher}), and it reads as follows.
\\
An orthogonal polynomial sequence, $\{B_n\}_{n\ge 0}$, is said
to be {\bf $q$-semiclassical} if
$$
\Phi(s)B^{[1]}_n(s) =\sum_{\nu=n-\sigma}^{n+t}\lambda_{n,\nu}
B_{\nu}(s), \ n\geq \sigma,\quad \lambda_{n,n-\sigma}\ne 0,
\ n\geq \sigma+1,
$$
where $\Phi$ is a polynomial of degree $t$ and $\sigma$ is a
non-negative integer such that $\sigma\ge \max\{t-2,0\}$.
\\
Recently, F. Marcell\'{a}n and R. Sfaxi \cite{SfMa2} have
established a second structure relation for the standard
semiclassical polynomials which reads as follows
\bt \label{the3.2p}
For any integer $\sigma\geq 0$, any monic polynomial $\Phi$, with
$\deg \Phi=t\leq \sigma+2$, and any {\rm SMOP} $\{B_n\}_{n\geq 0}$
with respect to a linear functional $u$, the following statements
are equivalent.
\be
\item[\bf(i)] There exist an integer $p\geq 1$ and an integer
$r\geq \sigma+t+1$, with $\sigma=\max(t-2,p-1)$, such that
$$
\sum_{\nu=n-\sigma}^{n+\sigma}\xi_{n,\nu}B_{\nu}(x)=
\sum_{\nu=n-t}^{n+\sigma}\varsigma_{n,\nu}B_{\nu}^{[1]}(x),\qquad
n\geq \max(\sigma,t+1),\eqno{(3.36)}
$$
where
$B_n^{[1]}(x)=(n+1)^{-1}B'_{n+1}(x)$,
\noindent\begin{align*}
\xi_{n,n+\sigma}=
\varsigma_{n,n+\sigma}&=\; 1,\;n\geq \max(\sigma,t+1),
\xi_{r,r-\sigma}\varsigma_{r,r-t}\ne 0, \\
\langle(\Phi u)',B_n\rangle &=\;0,\;\; p+1\leq
n\leq 2\sigma+t+1,\;\;
\langle(\Phi u)',B_p\rangle\ne 0,\;(\sigma\geq 1),
\end{align*}
and if $p=t-1$ then $\langle u,B_p^2\rangle^{-1}\langle u,\Phi B_p'
\rangle\notin \mathbb{N}^*$.
\item[\bf(ii)] The linear functional $u$ satisfies
$$
(\Phi u)'+\Psi u=0,
$$
where the pair $(\Phi,\Psi)$ is admissible, i.e.
the polynomial $\Phi$ is monic, $\deg \Phi=t$, $\deg \Psi=p\ge 1$
and if $p=t-1$ then $\frac 1{n!}\Psi^{(n)}(0)\not \in -\XX N^*$,
with associated integer $\sigma$.
\end{enumerate}
\et
\noindent
Now,
we are going to extend this result for the $q$-semiclassical
polynomials of the Hahn Tableau.
\\
Some years ago, P. Maroni and R. Sfaxi \cite{masf} introduced
the concept of diagonal sequence  for the standard semiclassical
polynomials.
The following definition extends this definition to the
$q$-semiclassical case.
\bd \label{d11}
Let $\{B_n\}_{n\geq 0}$ be a sequence of monic orthogonal
polynomials and $\phi$ a monic polynomial with $\deg\phi=t$.
When there exists an integer $\sigma\geq 0$ such that
\bq \label{defdiagseq}
\phi(s)B_n(s)=\sum_{\nu=n-\sigma}^{n+t}\theta_{n,\nu}
B_{\nu}^{[1]}(s),\; \; \theta_{n,n-\sigma}\ne 0,\quad
n\geq \sigma,
\eq
the sequence $\{B_n\}_{n\geq 0}$ is said to be
{\bf diagonal associated with $\phi$ and index $\sigma$.}
\ed \noindent
Obviously, the above finite-type relation, that we will call
diagonal relation, is nothing else that an example of second
structure relation for such a family. But, some
$q$-semiclassical orthogonal polynomials are not diagonal. As
an example, we can mention the case of a $q$-semiclassical
polynomial sequence $\{Q_n\}_{n\ge 0}$ orthogonal with
respect to the linear functional $v$, such that the
functional equation: $\dq v=\Psi v$, with
$\deg \Psi=2$, holds. In fact, the sequence $\{Q_n\}_{n\ge 0}$
satisfies the following relation
$$
(x(s+1)+v_{n,0})Q_n(s)=
qQ_{n+1}^{[1]}(s)+\rho_nQ_n^{[1]}(s),\ n\ge 0,
$$
where the lattice, $x(s)$, is $q$-linear, i.e.
$x(s+1)-qx(s)=\omega$,
$$
\ba{rl}
\rho_n= & \dfrac{q^{n+1}}{{\g C}}\dfrac{[n+1]}{\gamma_{n+1}},
\ n\ge 1,\quad \rho_0=0,\\[0.4cm] v_{n,0}=&\dfrac{\gamma_{n+2}
\gamma_{n+1}}{q^{n}[n+2]}{\g C}+\rho_n-q\beta_n-\omega,
\quad n\ge 0. \ea
$$
Here ${\g C}$ is a constant, $\gamma_n$ and $\beta_n$ are the
coefficients of the three-term recurrence relation (TTRR) that
the orthogonal polynomial sequence $\{Q_n\}_{n\ge 0}$
satisfies.
In fact, this sequence is not diagonal and it will be analyzed
more carefully in $\S$ \ref{4.1}.
\\
The aim of our contribution is to give, under certain
conditions, the second structure relation characterizing a
$q$-semiclassical polynomial sequence by a new relation
between the sequence of $q$-polynomials, $\{B_n\}_{n\geq 0}$,
and the polynomial sequence of monic $q$-differences,
$\{B^{[1]}_n\}_{n\geq 0}$, as follows
$$
\displaystyle\sum_{\nu=n-\sigma}^{n+\sigma}
\xi_{n,\nu}B_{\nu}(s)= \;\sum_{\nu=n-t}^{n+\sigma}
\varsigma_{n,\nu}B^{[1]}_n(s),\quad n\geq \max(t+1,\sigma),
$$
where $\xi_{n,n+\sigma}=\varsigma_{n,n+\sigma}=1,$ $n\geq
\max(t+1,\sigma)$, and there exists $r\geq \sigma+t+1$
such that $\xi_{r,r-\sigma} \varsigma_{r,r-t} \ne0$.
\\[0.2cm]
Notice that when $\sigma=0$ we get the second structure
relation \refe{sestre}.
\section{Preliminaries and notation} \noindent
Let $u$ be a linear functional in the linear space
$\mathbb{P}$ of polynomials with complex coefficients
and let $\mathbb{P}'$ be its algebraic dual space, i.e.,
the linear space of the linear functionals defined on
$\mathbb{P}$. We will denote by $\langle u,f\rangle$
the action of $u\in \mathbb{P}'$ on $f\in \mathbb{P}$
and by $(u)_{n}:=\langle u,x^n\rangle,\; n\geq 0$, the
moments of $u$ with respect to the sequence
$\{x^n \}_{n\geq 0}$.
\par\noindent Let us define the following operations in
$\mathbb{P}'$. For any polynomial $h$ and any
$c\in \mathbb{C}$, let $\dq u$, $hu$, and $(x-c)^{-1}u$
be the linear functionals defined on $\mathbb{P}$ by
(see \cite{mar6,kher})
\be
\item[\bf (i)] $\langle \dq u,f\rangle:=
-\langle u,\dq f\rangle, \quad f\in\mathbb{P}$,
\item[\bf (ii)]$\langle gu,f\rangle:=\langle u,gf\rangle,
\quad f,\;g\in\mathbb{P},$
\item[\bf (iii)] $\langle (x-c)^{-1}u,f\rangle:=
\langle u,\theta_c(f)\rangle, \quad
f\in\mathbb{P},\;c\in \mathbb{C},\quad \mbox{where}\;\;
\theta_c(f)(x)=\frac{f(x)-f(c)}{x-c}.$
\ee
Furthermore, for any linear functional $u$ and any polynomial $g$
we get
\bq \label{relgu}
L_{q,\omega}(gu):=\dq(gu)=g(q^{-1}(x-\omega))\dq u+
\dq(g(q^{-1}(x-\omega)))u.
\eq
Let $\{B_n\}_{n\geq 0}$ be a sequence of monic polynomials
(SMP) with $\deg B_n=n,\;n\geq 0$, and $\{u_n\}_{n\geq 0}$
its dual sequence, i.e. $u_n\in \mathbb{P}',\;n\geq 0$,
and $\langle u_n,B_m\rangle:=\delta_{n,m},\;n,
\;m\geq 0$, where $\delta_{n,m}$ is the Kronecker symbol.
The next results are very well-known
\cite{kher
}.
\bl \label{lem21} For any $u\in \mathbb{P}'$,
and any integer $m\geq 1$, the following statements are
equivalent.
\begin{enumerate}
\item[\bf (i)]  $\langle u,B_{m-1}\rangle\ne 0,\quad
\langle u, B_{n}\rangle=\;0,\;n\geq m.$
\item[\bf (ii)] There exist $\lambda_{\nu}\in \mathbb{C},$
$0\leq \nu\leq m-1,\; \lambda_{m-1}\ne 0,$ such that
$u=\dst \sum_{\nu=0}^{m-1}\lambda_{\nu}u_{\nu}.$
\end{enumerate}
\el \noindent
On the other hand, it is straightforward to prove \vspace{0.1cm}
\bl \label{lem22}
For any $(\widehat t,\widehat \sigma,\widehat r)\in \XX N^3$,
$\widehat r\ge \widehat \sigma+ \widehat t+1$ and any
sequence of monic polynomials $\{\Omega_n\}_{n\ge 0}$,
$\deg \Omega_n=n$, $n\ge 0$, with dual sequence
$\{w_n\}_{n\ge 0}$ such that
$$\ba{rl}
\Omega_n(x)=& \dst \sum_{\nu=n-\widehat t}^n
\lambda_{n,\nu}B_\nu(x), \quad n\ge \widehat t+\widehat
\sigma+1,\quad \lambda_{\widehat r,\widehat r-\widehat t}
\ne 0, \\[0.5cm]  \Omega_n(x)=&B_n(x),\quad 0\le n\le
\widehat t+\widehat \sigma,
\ea
$$
we have that $w_k=u_k$ for every $0\le k\le \widehat \sigma$.
\el \noindent
The linear functional $u$ is said to be quasi-definite if, for
every non-negative integer, the leading principal Hankel
submatrices $H_n\!=\!{\big(}(u)_{i+j}{\big)}_{i,j=0}^{n}$
are non-singular for every $n\ge 0$.
Assuming $u$ is quasi-definite, there exists a sequence of monic
polynomials $\{B_n\}_{n\geq 0}$ such that (see \cite{chi})
\begin{enumerate}
\item[\bf (i)] $\deg B_n=n,\;n\geq 0$,
\item [\bf (ii)] $\langle u,B_nB_m\rangle=r_n\delta_{n,m}$,
with $r_n=\;\langle u,B_n^2\rangle\ne 0,\; n\geq0$.
\end{enumerate}
The sequence $\{B_n\}_{n\geq 0}$ is said to be the sequence of
monic orthogonal polynomials, in short SMOP with respect to
the linear functional $u$.
\\
If $\{B_n\}_{n\geq 0}$ is a SMOP, with respect to the
quasi-definite linear functional $u$, then it is well-known
(see \cite{mar6}) that its corresponding dual sequence
$\{u_n\}_{n\geq 0},$ is
\bq \label{2.1}
u_n=\;r_n^{-1}B_n u,\; n\geq 0.
\eq
\bn
We assume $u_0=u$, i.e. the linear functional $u$ is normalized.
\en
On the other hand, (see \cite{chi}), the sequence
$\{B_n\}_{n\geq 0}$ satisfies a three-term recurrence
relation (TTRR)
\bq \label{2.2}
B_{n+1}(x)=(x-\beta_n)B_n(x)-\gamma_nB_{n-1}(x),\; n\geq
0,
\eq
with $\gamma_n\ne 0,\; n\geq 1$ and
$B_{-1}(x)=0,\; B_{0}(x)=1$.
\\
Conversely, given a SMP,
$\{B_n\}_{n\geq 0}$, generated by a
recurrence relation \refe{2.2} as above with $\gamma_n\ne 0,
\; n\geq 1,$ there exists a unique normalized quasi-definite linear
functional $u$ such that the family $\{B_n\}_{n\geq 0}$ is
the corresponding SMOP. This result is known as Favard
Theorem (see \cite{chi}).
\\
An important family of linear functionals is constituted by
the $q$-semiclassical linear func\-ti\-o\-nals, i.e., when
$u$ is quasi-definite and satisfies 
\bq \label{qdisteq}
\dq(\Phi u)=\Psi u.
\eq
Here $(\Phi,\Psi)$
is an admissible pair of polynomials, i.e., the
polynomial $\Phi$ is monic, $\deg \Phi=t$,
$\deg \Psi=p\geq 1$, and if  $p=t-1,$ then the following
condition 
holds
$$
\lim_{q\uparrow 1}\frac 1 {[p]!}[\dq]^p\,\Psi(0):=
\lim_{q\uparrow 1} \frac 1 {[p]!}\overbrace{\Delta^{(1)}
\cdots \Delta^{(1)}}^{p} \Psi(0)\ne -n, \quad n\in \XX N^*,
$$
where $[m]!=[1][2]\cdots[m]$, $m\in \XX N^*$, is the $q$-analog
of the usual factorial.\\
\noindent The pair $(\Phi,\Psi)$ is not unique. In fact, under
certain conditions \refe{qdisteq} can be simplified, so we
define the class of $u$ as the minimum value of
$\max\big( \deg(\Phi)-2,\;\deg(\Psi)-1\big),$ for all
admissible pairs $(\Phi,\Psi).$ The pair $(\Phi,\Psi)$ giving
the class $\sigma$ $(\sigma\geq0$ because $\deg(\Psi)\geq 1)$
is unique \cite{kher}.
\par\noindent When $u$ is $q$-semiclassical of
class $\sigma,$ the corresponding SMOP is said to be
$q$-semiclassical of class $\sigma$.
\\
When $\sigma=0,$ i.e., $\deg \Phi\leq 2$ and $\deg \Psi=1,$
then $u$ is $q$-classical (Askey-Wilson, $q$-Racah, Big
$q$-Jacobi, $q$-Charlier, etc). For more details see
\cite{ks,medrenmar,nsu}.
\section{Main results} \noindent
First, we will present particular cases of {\it diagonal sequences}.
\\
Let $\{P_n\}_{n\ge 0}$ and $\{Q_n\}_{n\ge 0}$ be sequences of
monic polynomials, $\{v_n\}_{n\ge 0}$ and $\{w_n\}_{n\ge 0}$
their corresponding dual sequences. Let $\phi$ be a
monic polynomial of degree $t$.
\bd
The sequence $\{P_n\}_{n\ge 0}$ is said to be compatible with
$\phi$ if $\phi v_n\ne 0$, $n\ge 0$.
\ed
\bl \label{prop21}\cite[Prop. 2.1]{mar6} Let $\phi$ be as above.
For any sequence $\{P_n\}_{n\ge 0}$ compatible with $\phi$, the
following statements are equivalent.
\be
\item[\bf (i)] There is an integer $\sigma\ge 0$ such that
\begin{eqnarray}
\label{eq1} \dst \phi(x)Q_n(x)=\sum_{\nu=n-\sigma}^{n+t}
\lambda_{n,\nu} P_{\nu}(x),\quad n \ge \sigma,\\[0.3cm]
\label{eq2} \exists \, r\ge \sigma\ : \quad
\lambda_{r,r-\sigma}\ne 0.
\end{eqnarray}
\item[\bf (ii)] There are an integer $\sigma\ge 0$ and a
mapping from $\XX N$ into $\XX N:m\mapsto \mu(m)$ satisfying
\begin{eqnarray}
\label{eq3} \max\{0,m-t\}\le \mu(m)\le m+\sigma,
\quad m\ge 0, \\[0.3cm] \label{eq4}\exists \
m_0\ge 0\quad {\rm with} \quad \mu(m_0)=m_0+\sigma,
\end{eqnarray}
such that
\bq \label{eq5}
\ba{l}\dst \phi v_m=\sum_{\nu=m-t}^{\mu(m)}
\lambda_{\nu,m} w_\nu, \quad m\ge t,\\[0.4cm]
\lambda_{\mu(m),m}\ne 0,\quad m\ge 0.
\ea
\eq
\ee
\el
\bp \cite[Prop. 2.2]{mar6} Assume $\{Q_n\}_{n\ge 0}$ is
orthogonal and $\{P_n\}_{n\ge 0}$ is compatible with $\phi$.
Then the sequences $\{P_n\}_{n\ge 0}$ and $\{Q_n\}_{n\ge 0}$
fulfil the finite-type relations \refe{eq1}-\refe{eq2} if
and only if there are an integer $\sigma\ge 0$ and a mapping
from $\XX N$ into $\XX N:m\mapsto \mu(m)$ satisfying
\refe{eq3} and \refe{eq4}.
Moreover, there exist $\{k_m\}_{m\ge 0}$ and a sequence
$\{\Lambda_{\mu(m)}\}_{m\ge 0}$ of monic polynomials with
$\deg(\Lambda_{\mu(m)})=\mu(m)$, $m\ge 0$, such that
\bq \label{relomega}
\phi v_m=k_m\Lambda_{\mu(m)} w_0,\quad m\ge 0.
\eq
\ep
\noindent From these two results we get
\bc \label{lem16} \cite[Prop. 1.6]{masf} Let $\phi$ be
as above.
For sequences of monic orthogonal polynomials (SMOP)
$\{P_n\}_{n\ge 0}$ and $\{B_n\}_{n\ge 0}$ orthogonal
with respect to linear functionals $v$ and $u$,
respectively, the following statements are equivalent.
\be
\item[\bf (i)] There exists an integer $\sigma\ge 0$ such that
$$
\phi(s)P_n(s)=\sum_{\nu=n-\sigma}^{n+t}\lambda_{n,\nu}
B_\nu^{[1]}(s),\quad \lambda_{n,n-\sigma}\ne 0, \ n\ge \sigma.
$$
\item[\bf (ii)] There exists a monic polynomial sequence
$\{\Omega_{n+\sigma}\}_{n\ge 0}$, with
$\deg(\Omega_{n+\sigma})=n+\sigma$, $n\ge 0$ and non-zero
constants $k_n$, $n\ge 0$ such that
\bq \label{relmom}
\phi u_n^{[1]}=k_n \Omega_{n+\sigma} v_0.
\eq
where $\{u_n^{[1]}\}_{n\ge 0}$ is the dual sequence of
$\{B_n^{[1]}\}_{n\ge 0}$.
\ee
\ec \noindent
Thus we can prove
\bp \label{diag_seq} Any diagonal sequence, $\{B_n\}_{n\ge0}$,
orthogonal with respect a linear functional $u$ is necessarily
semiclassical and $u$ satisfies
\bq \label{diaeq}
\dq (\phi(qx+\omega)\Omega_{n+\sigma}(x)u)=
\psi_n(x)u,\ n\ge 0,
\eq
where
\bq \label{psi}
\psi_n(s)=\frac{\phi(s+1)-\phi(s-1)}{\Delta x(s)}\,
\Omega_{n+\sigma}(s)-d_n \phi(s) \phi(s-1)B_{n+1}(s),
\eq
and
\bq \label{kn}
d_n=[n+1]\frac{\pe {u}{B_{n+\sigma}^2}}{
\pe {u}{B_{n+1}^2} \lambda_{n+\sigma,n}},\ n\ge 0.
\eq
Furthermore, the sequence $\{\Omega_{n+s}\}_{n\ge 0}$ satisfies
\bq \label{Omerel}
\Omega_{n+\sigma}(s)\dq \Omega_{\sigma}(s)-
\Omega_{\sigma}(s)\dq\Omega_{n+\sigma}(s)=
\phi(s+1)\{d_n \Omega_\sigma(s)B_{n+1}(s+1)-
d_0\Omega_{n+\sigma}(s)B_1(s+1)\}.
\eq
\ep
\bdm
Let $\{B_n\}_{n\ge 0}$ be a diagonal sequence in the sense of
Definition \ref{d11} and assume the linear
functional $u$ is normalized. Then from Lemma
\ref{lem16} there exist a sequence of monic polynomials
$\{\Omega_{n+\sigma}\}_{n\ge 0}$ and non-zero constants
$\{k_n\}_{n\ge 0}$ such that
$$
\phi u_n^{[1]}=k_n\Omega_{n+\sigma} u.
$$
Then
\bq \label{eqa1}\begin{split}
k_n\dq(\Omega_{n+\sigma} u)=\ & \dq(\phi(q^{-1}(x-
\omega)))u_n^{[1]}+\phi(q^{-1}(x-\omega))\dq u_n^{[1]}\\
=\ &\dq(\phi(q^{-1}(x-\omega)))u_n^{[1]}-\frac{[n+1]}
{\langle u,B_{n+1}^2 \rangle}\,\phi(q^{-1}(x-\omega))
B_{n+1}(x)u(s), \end{split}
\eq
as well as
\bq \label{relPhi}
\dq\left(\phi(s)\phi(s-1)\right)=\phi(s)
\frac{\phi(s+1)-\phi(s-1)}{\Delta x(s)}.
\eq
Combining \refe{eqa1} and \refe{relPhi}, a straightforward
calculation yields \refe{diaeq}, \refe{psi}, and \refe{kn}.
\\
Taking \refe{diaeq} for $n=0$ and cancelling out
$\dq (\phi(qx+\omega)u)$, from the quasi-definite character
of $u$ we obtain \refe{Omerel}.
\edm
\bc \label{cor1}\cite[Corollary 2.3]{masf} If
$\{B_n\}_{n\ge 0}$ is a diagonal sequence given by
\refe{defdiagseq}, then we get
\bq \label{relpar}
\half t\le \sigma\le t+2.
\eq
\ec \noindent
For a linear functional $u$, let $(\Phi,\Psi)$ be the
minimal admissible pair of polynomials with $\Phi$ monic,
$\deg \Phi=t$, and $\deg \Psi=p\geq 1$, defined as above.
To this pair we can associate the non-negative integer
$\sigma:=\max(t-2,p-1)\geq 0$.
\\
Now, given  $\{B_n\}_{n\geq 0}$, a SMOP with respect
to $u$, we get
\bq \label{eqsem}
\Phi(s) B_{n}^{[1]}(s)=\;\displaystyle
\sum_{\nu=0}^{n+t}\lambda_{n,\nu}B_{\nu}(s),
\quad n\geq\max(t-1,0),
\eq
where $\lambda_{n,n+t}=1$ and
$$\ba{rl}
\lambda_{n,\nu}  = & \!\!\!\!r_\nu^{-1}\langle u,
\Phi(s) B_{n}^{[1]}(s)B_\nu(s)\rangle= \dfrac {r_\nu^{-1}}{[n+1]}
\langle B_{\nu}\Phi u, \dq B_{n+1}\rangle\\[0.3cm]
= & \!\!\!\!-\dfrac {r_\nu^{-1}}{[n+1]} \langle B_\nu(q^{-1}
(x-\omega)) \dq(\Phi u) \!+\!\dq(B_\nu(q^{-1}
(x-\omega)))\Phi u,B_{n+1}\rangle, \ 0\leq\nu\leq n+t.
\ea
$$
\bl \label{lem31}\cite[Prop. 3.2]{kher} For any
monic polynomial $\Phi$, $\deg \Phi=t$,
and any {\rm SMOP} $\{B_n\}_{n\geq 0}$ with respect to
$u$, the following statements are equivalent.
\renewcommand\theenumi{(\roman{enumi})}
\renewcommand\labelenumi{\theenumi}
\be
\item
There exists a non-negative integer $\sigma$ such that
\bq \label{cond1}
\Phi(s) B_{n}^{[1]}(s)=\;\displaystyle
\sum_{\nu=n-\sigma}^{n+t}\lambda_{n,\nu}B_{\nu}(s),
\quad n\geq \sigma,
\eq
\bq \label{cond2}
\lambda_{n,n-\sigma}\ne 0,\quad n\geq \sigma+1.
\eq
\item
There exists a polynomial $\Psi$, $\deg \Psi=p\geq 1$,
such that
\bq \label{3.8}
\dq(\Phi u)=\Psi u.
\eq
where the pair $(\Phi,\Psi)$ is admissible.
\item
There exist a non-negative integer $\sigma$ and a polynomial
$\Psi$, with $\deg \Psi=p\geq 1$, such that
\bq \label{3.9}
\Phi(s)\dq B_n(s-1)+\Psi(s)B_n(s-1)=
\displaystyle\sum_{\nu=n-t}^{n+\sigma(n)}
{\tilde\lambda}_{n,\nu}B_{\nu+1}(s),\quad n\geq t,
\eq
\bq \label{lambtilden}
{\tilde\lambda}_{n,n-t}\ne 0,\quad n\geq t,
\eq
where $\sigma=\max(p-1,t-2)$, the pair $(\Phi,\Psi)$
is admissible, and
\bq \label{lambtilde0}
\sigma(n)=\left\{\ba{rl}p-1,& n=0, \\ \sigma, & n\ge 1.
\ea\right.
\eq
\ee
We can write
\bq \label{latim}
\tilde \lambda_{n,\nu}=-[\nu+1]\frac{\langle u,
B_n^2\rangle} {\langle u,B_{\nu+1}^2\rangle}
\lambda_{\nu,n}, \quad 0\leq \nu\leq n+\sigma.
\eq
\el
\bdm
(i)$\Rightarrow$ (ii), (iii). Assuming (i), from Lemma
\ref{prop21} and taking $P_n=B_n$ and $Q_n=B_n^{[1]}$,
we get
$$
\Phi u_m=\sum_{\nu=0}^{\mu(m)}\lambda_{\nu,m}u_\nu^{[1]},
\ m\ge 0.
$$
On the other hand, \refe{cond2} implies $\mu(m)=m+\sigma$, $m\ge 1$.
\\
Taking into account that
\bq \label{cond3}
\dq u_m^{[1]}=-[m+1]u_{m+1},\ m\ge 0,
\eq
we have
$$
\dq (\Phi u_m)=-\sum_{\nu=0}^{\mu(m)}\lambda_{\nu,m}
[\nu+1]u_{\nu+1},\ m\ge 0.
$$
In accordance with the orthogonality of $\{B_n\}_{n\ge 0}$,
we get
\bq \label{aux1}
\dq(\Phi B_m u)=-\Psi_{\mu(m)+1} u,
\ m\ge 0,
\eq
with
\bq \label{Psi}
\Psi_{\mu(m)+1}(s)=\sum_{\nu=0}^{\mu(m)}\lambda_{\nu,m}
[\nu+1]B_{\nu+1}(s),\ m\ge 0.
\eq
Taking  $m=0$ in \refe{aux1}, we have
\bq \label{aux2}
\dq(\Phi u)=-\Psi_{\mu(0)+1} u.
\eq
Inserting \refe{aux2} in \refe{aux1} and because $u$
is quasi-definite, we get
$$
\Phi(s)\dq B_m(s-1)-\Psi_{\mu(0)+1}(s)B_m(s-1)=
-\Psi_{\mu(m)+1}(s),\ m\ge 0.
$$
The consideration of the degrees in both hand sides leads to
\bi
\item If $t-1>\mu(0)+1$, which implies $t\ge 3$, then
$t=\sigma+2$, $\mu(0)<\sigma$.
\item If $t-1\le \mu(0)+1$, then $\mu(0)=\sigma$,
$t\le \sigma+2$.
\ei
Obviously, the pair $(\Phi,-\Psi_{\mu(0)+1})$ is
admissible and putting $p=\mu(0)+1$, we have
$\sigma=\max(p-1,t-2)$. So \refe{3.9} and
\refe{lambtilden} are valid from \refe{latim}.
\\
Thus, we have proved that (i)$\Rightarrow$(ii) and
(i)$\Rightarrow$(iii).
\\
(ii)$\Rightarrow$(iii). Consider $m\ge 0$. Thus
$$
\Phi(s)\dq B_m(s-1)+\Psi(s)B_m(s-1)=
\dst\sum_{\nu=0}^{m+\sigma(m)+1} \lambda'_{m,\nu}B_{\nu}(s).
$$
We successively derive from this
$$
\pe {u} {(\Phi(s)\dq B_m(s-1)+\Psi(s)B_m(s-1))B_\mu}=
\lambda'_{m,\mu}\pe {u} {B^2_\mu},\ 0\le \mu\le m+\sigma+1.
$$
A straightforward calculation yields
\bq \label{aux6}
\pe {u} {(\Phi(s)\dq B_m(s-1)+\Psi(s)B_m(s-1))B_\mu}=
-\pe {u} {\Phi(s)B_m(s)\dq B_\mu(s)}.
\eq
Then
$$
-\pe {u} {\Phi(s)B_m(s)\dq B_\mu(s)}=
\lambda'_{m,\mu}\pe {u} {B^2_\mu}.
$$
Consequently, $\lambda'_{m,\mu}=0$, $0\le \mu\le m-t$,
$\lambda'_{m,0}=0$, $m\ge 0$. Moreover,
for $\mu=m-t+1$, $m\ge t$,
$$
-\pe {u} {\Phi(s)P_m(s)\dq P_{m-t+1}(s)}=
-[m-t+1]\pe {u} {B_m^2}=
\lambda'_{m,m-t+1}\pe {u} {B^2_{m-t+1}}.
$$
Therefore, for $m\ge t$,
$$
\Phi(s)\dq B_m(s-1)+\Psi(s)B_m(s-1)=\dst\sum_{\nu=m-t}^{m+
\sigma(m)} \lambda'_{m,\nu+1}B_{\nu+1}(s),\
\lambda'_{m,m-t+1}\ne 0.
$$
(iii)$\Rightarrow$(i). From \refe{3.9}, we get
$$\ba{rl}
\dst \sum_{\nu=0}^{m+\sigma(m)} {\tilde\lambda}_{m,\nu}
\delta_{n,\nu+1} = &  \pe {u_n} {\Phi(s)\dq B_m(s-1)+
\Psi(s)B_m(s-1)} \\
= & -\pe {\dq(\Phi u_n)-\Psi u_n} {B_m(s-1)}.
\ea
$$
For $n=0$, $\pe {\Psi u-\dq(\Phi u)}
{B_m(s-1)}=0$, $m\ge 0$. Therefore
\bq \label{aux3}
\dq(\Phi u)=\Psi u.
\eq
Moreover, using \refe{aux6} and the orthogonality of
$\{B_n\}_{n\ge 0}$, we get
$$
\pe {u_n} {\Phi(s)\dq B_m(s-1)+\Psi(s)B_m(s-1)}=
 -r_n^{-1}\pe {u} {\Phi(s)B_m(s)\dq B_n(s)}.
$$
Furthermore, making $n\to n+1$, we obtain
$$
\left\{\ba{rl}
\pe {(\Phi\dq B_{n+1})u} {B_m}&\hspace{-0.3cm}=
0, \ m\ge n+t+1,\ n\ge 0, \\[0.3cm] \pe {(\Phi\dq
B_{n+1})u} {B_{n+t}} & \hspace{-0.3cm}= -r_{n+1}
\tilde \lambda_{n+t,n}\ne 0,\  n\ge 0.
\ea\right.
$$
According to Lemma \ref{lem21},
$$
(\Phi\dq B_{n+1})u=
-\sum_{\nu=n-\sigma}^{n+t} r_n\tilde \lambda_{\nu,n}u_\nu,\
n\ge \sigma.
$$
The orthogonality of $\{B_n\}_{n\ge 0}$ leads to
$$
(\Phi\dq B_{n+1})u=
-\!\!\! \sum_{\nu=n-\sigma}^{n+t} \!\!\!\left(\tilde
\lambda_{\nu,n} \frac{\pe {u}{B_{n+1}^2}}
{\pe {u}{B_{\nu}^2}}B_\nu \right)u ,\ n\ge 0.
$$
From \refe{aux3} and taking into account $u$ is
quasi-definite, we finally obtain
\refe{cond1}--\refe{cond2} in accordance with \refe{latim}.
\edm
\\
In an analog way we can prove the following result
\bl \label{lem31p}\cite[Lemma 3.1]{SfMa2} For any monic
polynomial $\Phi$, $\deg \Phi=t$, and any {\rm SMOP}
$\{B_n\}_{n\geq 0}$ with respect to $u$, the
following statements are equivalent.
\renewcommand\theenumi{(\roman{enumi})}
\renewcommand\labelenumi{\theenumi}
\be
\item
There exists a non-negative integer $\sigma$ such that
the polynomials $B_n$ satisfy
\bq \label{cond1p}
\dq(\Phi(s-1) B_{n}(s))=\;\displaystyle
\sum_{\nu=n-\sigma-1}^{n+t-1}\lambda_{n,\nu}
B_{\nu}(s),\quad n\geq \sigma+1,
\eq
\bq \label{cond2p}
\lambda_{n,n-\sigma-1}\ne 0,\quad n\geq t+\sigma+2.
\eq
\item
There exists a polynomial $\Psi$, $\deg \Psi=p\geq 1$,
such that
\bq \label{3.8p}
\dq(\Phi u)=\Psi u.
\eq
where the pair $(\Phi,\Psi)$ is admissible.
\item
There exist a non-negative integer $\sigma$ and a
polynomial $\Psi$, $\deg \Psi=p\geq 1$, such that
\bq \label{3.9p}
\Phi(s)\dq B_n(s-1)+\Psi(s)B_n(s-1)-B_n(s)\dq \Phi(s-1)=
\displaystyle\sum_{\nu=n-t+1}^{n+\sigma(n)+1}
{\tilde\lambda}_{n,\nu}B_{\nu}(s),\quad n\geq t,
\eq
\bq \label{lambtildenp}
{\tilde\lambda}_{n,n-t+1}\ne 0,\quad n\geq t,
\eq
where $\sigma=\max(p-1,t-2)$ and the pair $(\Phi,\Psi)$
is admissible.
We can write
\bq \label{latimp}
\tilde \lambda_{n,\nu}=-\frac{\langle u,
B_m^2\rangle} {\langle u,B_{\nu}^2\rangle}
\lambda_{\nu,n}, \quad 0\leq \nu\leq n+\sigma(n)+1, \ n\ge 0.
\eq
\ee
\el
\subsection{First Characterization of
$q$-semiclassical polynomials}
\bt \label{theo31}For a monic polynomial $\Phi$, $\deg \Phi=t$,
and any SMOP $\{B_n\}_{n\ge 0}$ with respect to $u$, the
following statements are equivalent.
\be
\item[(i)] There exist a non-negative integer $\sigma$, an
integer $p\ge 1$, and an integer $r\ge \sigma+t+1$, with
$\sigma=\max(t-2,p-1)$, such that
\bq \label{relteo1}
\sum_{\nu=n-\sigma}^{n+t} \alpha_{n,\nu}B_\nu(s)=
\sum_{\nu=n-t}^{n+t} v_{n,\nu}B^{[1]}_\nu(s), \quad n\ge
\max(\sigma,t),
\eq
where $\alpha_{n,n+t}=v_{n,n+t}=1$, $n\ge \max(\sigma,t)$,
$\alpha_{r,r-\sigma}v_{r,r-t}\ne 0$,
$$
\langle\dq (\Phi u),B_n\rangle=0, \ p+1\le n\le \sigma+
2t+1, \quad \langle \dq(\Phi u),B_p\rangle\ne 0,
$$
and if $p=t-1$, then $\dst \lim_{q\uparrow 1}
\langle u,B_p^2\rangle^{-1}
\langle u,\Phi \dq B_p\rangle\ne -m$, $m\in \XX N^*$.
\item[(ii)]There exists a polynomial $\Psi$,
$\deg \Psi=p\ge 1$, such that
$$
\dq (\Phi u)=\Psi u,
$$
and the pair $(\Phi,\Psi)$ is admissible.
\ee
\et
\bdm
$(i)\Rightarrow (ii)$. Consider the SMP
$\{\Omega_{n}\}_{n\ge 0}$ defined by
$$
\ba{rl}
\Omega_{n+t+1}(s)=& \dst \sum_{\nu=n-t}^{n+t} \frac{[n+t+1]}
{[\nu+1]}v_{n,\nu} B_{\nu+1}(s),\quad n\ge \sigma+t+1,
\\[0.4cm] \Omega_n(s)= & B_n(s), \quad 0\le n\le \sigma+2t+1.
\ea
$$
From \refe{relteo1},
\bq \label{au1}
\dq(\Omega_{n+t+1}(s))=[n+t+1]\sum_{\nu=n-\sigma}^{n+t}
\alpha_{n,\nu}B_\nu(s),\quad n\ge \sigma+t+1.
\eq
Since $u$ is quasi-definite, then
$$
\ba{rl}
\langle \dq(\Phi u),\Omega_{n+t+1}\rangle & = -\langle u,
\Phi \dq \Omega_{n+t+1}\rangle \\
& = \dst -[n+t+1]\sum_{\nu=n-\sigma}^{n+t} \alpha_{n,\nu}
\langle u,\Phi B_\nu\rangle=0, \quad n\ge \sigma+t+1.
\ea
$$
Therefore, $\pe {\dq(\Phi u)} {\Omega_{n}}=0$, $n\ge \sigma+2t+1$,
 and by hypothesis $\langle \dq(\Phi u),\Omega_{n}
\rangle=0$, $p+1\le n\le \sigma+2t+1$, then
$\pe {\dq(\Phi u)} {\Omega_{n}}=0$ for $n\ge p+1$, and
$\langle \dq(\Phi u), \Omega_{p}\rangle\ne 0$. Hence, if we
denote $\{w_n\}_{n \ge 0}$ the dual sequence of
$\{\Omega_n\}_{n\ge 0}$ and apply Lemma \ref{lem21},
then
\bq \label{au2}
\dq(\Phi u)=\sum_{\nu=1}^p \langle \dq(\Phi u),B_\nu\rangle
w_\nu.
\eq
On the other hand, if we take $\widehat t=2t$,
$\widehat \sigma=\sigma+1$, and $\widehat r=r+t+1$, then
$$
\ba{rl}\Omega_{n}(s)=& \dst \sum_{\nu=n-\widehat t}^n\tilde
v_{n,\nu}B_\nu(s), \quad n\ge \widehat \sigma+\widehat t+1,
\\[0.4cm] \Omega_n(s)= & B_n(s), \quad 0\le n\le \widehat
\sigma+\widehat t,\ea
$$
where
$$
\ba{rl} \dst \widetilde v_{n,\nu}=\frac{[n]}{[\nu]}\,
v_{n-t-1,\nu-1}, \quad n-\widehat t\le \nu\le n,
\quad n\ge \widehat \sigma+\widehat t+1,\\[0.4cm]
\widetilde v_{\widehat r,\widehat r-\widehat t}=\dst
\frac{[r+t+1]}{[r-t+1]}\,v_{r,r-t}\ne 0, \quad \widehat r \ge
\sigma + 2t+2=\widehat \sigma+\widehat t+1.
\ea
$$
From Lemma \ref{lem22} and \refe{2.1}, it follows  that
$w_k=u_k=\langle u,B_k^2\rangle^{-1}B_k$, $0\le k\le \widehat
\sigma=\sigma+1$. So, relation \refe{au2} becomes
$$
\dq(\Phi u)=\Psi u,
$$
where
$$
\Psi(s)=-\sum_{\nu=1}^p \langle u,B_\nu^2
\rangle^{-1} \langle u,\Phi\dq B_\nu\rangle B_\nu(s),
$$
with $\deg \Psi=p$, as well as we have
$\langle u,\Phi \dq B_p\rangle\ne 0$ and, as a consequence, the pair
$(\Phi,\Psi)$ is admissible with associated integer $\sigma$.
\\
$(ii)\Rightarrow (i)$. From Lemma \ref{lem31p} (i) and making
$n\to n+1$ we have
\bq\label{auu1}
\dq(\Phi(s-1) B_{n+1}(s))=\sum_{\nu=n-\sigma}^{n+t}
\lambda_{n+1,\nu}B_\nu(s), \quad n\ge \sigma,
\eq
where $\lambda_{n+1,n+t}=[n+t+1]$, $n\ge \sigma$, and
$\lambda_{n+1,n-\sigma}\ne 0$, $n\ge t+\sigma+1$.
\\
On the other hand, the orthogonality of $\{B_n\}_{n\ge 0}$
yields
$$
\Phi(s-1)B_{n+1}(s)=\sum_{\nu=n-t}^{n+t}\frac{\langle u,
\Phi(s-1) B_{n+1}(s)B_{\nu+1}(s)\rangle}{\langle u,
B^2_{\nu+1}\rangle}B_{\nu+1}(s),\quad n\ge t-1.
$$
Hence,
\bq \label{auu2}
\dq(\Phi(s-1)B_{n+1}(s))=\sum_{\nu=n-t}^{n+t}\frac{[\nu+1]
\langle u, \Phi(s-1) B_{n+1}(s)B_{\nu+1}(s)\rangle}
{\langle u, B^2_{\nu+1} \rangle}B^{[1]}_\nu(s),\quad n\ge t.
\eq
From  \refe{auu1} and \refe{auu2}, we obtain \refe{relteo1}
with
$$
\ba{rl}
\alpha_{n,\nu}=& \dst \frac{\lambda_{n+1,\nu}}{[n+t+1]},
\quad n-\sigma\le \nu\le n+t, \\[0.3cm]
v_{n,\nu}=& \dst \frac{[\nu+1]\langle u,\Phi(s-1)
B_{n+1}(s)B_{\nu+1}(s)\rangle}{[n+t+1] \pe u {B_{\nu+1}^2}},
\quad n-t\le \nu\le n+t,\\[0.4cm] \alpha_{n,n-\sigma}
v_{n,n-t}\ne & \dst 0,\quad n \ge \sigma + t+1.
\ea
$$
Then,
$$
\langle \dq(\Phi u),B_n\rangle = -\langle u,\Phi \dq
B_n\rangle = \left\{\ba{l}0, \ p+1\le n\le \sigma+2t+1,
\\[0.3cm] \dst  \frac 1{[p]!} [\Delta^{(1)}]^p
\Psi(0)\langle u,B_p^2
\rangle, \ n=p=\deg \Psi,\ea \right.
$$
and if $p=t-1$, the $q$-admissibility of $(\Phi,\Psi)$ yields
$\dst \lim_{q\uparrow 1}\langle u,B_p^2\rangle^{-1}\langle
u,\Phi \dq B_p\rangle\ne -m, \ m\in \XX N^*$.
\hspace{-0.1cm}
\edm
\\
In the case of $q$-classical linear functionals, we get
the following result
\bc Let $\{B_n\}_{n\ge 0}$ be a SMOP with respect to $u$,
and a monic polynomial $\Phi$, with $\deg \Phi=t\le2$,
such that $\pe u \Phi\ne 0$, then the following statements
are equivalent.
\be
\item[\bf (i)] The linear functional $u$ is $q$-classical, i.e.
there exists a polynomial $\Psi$ with $\deg \Psi=1$ such
that $\dq(\Phi u)=\Psi u$.
\item[\bf (ii)]$\sum_{\nu=n}^{n+t}
\alpha_{n,\nu}B_\nu(s)= \sum_{\nu=n-t}^{n+t}v_{n,\nu}
B^{[1]}_\nu(s)$, $n\ge t$. Furthermore, there exists
an integer $r\ge t+1$ such that $\alpha_{r,r}v_{r,r-t}\ne 0$,
and if $t=2$ then $\lim_{q\uparrow 1}\pe u {B_1^2}^{-1}
\pe u {\Phi}\ne -m$, $m\in \mathbb{N}^*$.
\ee
\ec
\subsection{Second Characterization of $q$-semiclassical
polynomials}
From the previous characterization, we can not recover the
second structure relation of $q$-classical orthogonal
polynomials \refe{sestre}.
Our goal is to establish the characterization that allows
us to deduce such a case.
\\
First, we have the following result.
\bp \label{prop34}
For any monic polynomial $\Phi$, with $\deg \Phi=t$, and any
SMOP $\{B_n\}_{n\ge 0}$ with respect to $u$, the following
statements are equivalent.
\be
\item[\bf (i)] There exists  a polynomial $\Psi$,
$\deg \Psi=p\ge 1$, such that
\bq \label{3.30}
\dq (\Phi u)=\Psi u,
\eq
where the pair $(\Phi,\Psi)$ is admissible.
\item[\bf (ii)] There exist a non-negative integer $\sigma$
and a polynomial $\Psi$, with $\deg \Psi=p\ge 1$, such that
\bq \label{3.31}
\Phi(s)[\dq]^2 B_n(s-1)+\dq(\Psi(s)B_n(s-1))-B_n(s)[\dq]^2
\Phi(s-1)= \dst\sum_{\nu=n-\sigma}^{n+\sigma(n)}
{\vartheta}_{n,\nu}B_{\nu}(s),\ n\geq \sigma,
\eq
where
$\vartheta_{n,n-\sigma} \ne 0$ either
$n\ge \sigma+t+1$ or $n=\sigma+t$ and $p\ge t-1$,\ 
$\sigma=\max(t-2,p-1)$, and the pair $(\Phi,\Psi)$ is
admissible. We can write
\bq
\vartheta_{n,\nu}=\frac{\pe u {B_n^2}}{\pe u {B_\nu^2}}\,
\vartheta_{\nu,n}, \quad 0\le \nu\le n+\sigma(n),\ n\ge 0.
\eq
\ee
\ep
\bdm
We have
\bq \label{3.35}
\Phi(s)[\dq]^2 B_n(s-1)+\dq(\Psi(s)B_n(s-1))-B_n(s)[\dq]^2
\Phi(s-1)= \dst\sum_{\nu=0}^{n+\sigma(n)}
{\vartheta}_{n,\nu}B_{\nu}(s),\ n\geq 0,
\eq
where for all integers $0\le \nu\le n+\sigma(n)$, and $n\ge 0$,
$$
\pe u {B_\nu^2} \vartheta_{n,\nu}=\pe u
{(\Phi(s)[\dq]^2 B_n(s-1)+\dq(\Psi(s)B_n(s-1))-
B_n(s)[\dq]^2 \Phi(s-1))B_\nu}.
$$
Taking into account \refe{2.1} and \refe{3.31}, a straightforward
calculation leads to
$$
\pe u {B_\nu^2} \vartheta_{n,\nu}=\pe u
{(\Phi(s)[\dq]^2 B_\nu(s-1)+\dq(\Psi(s)B_\nu(s-1))-
B_\nu(s)[\dq]^2 \Phi(s-1))B_n}.
$$
Therefore, inserting \refe{3.35}
$$
\pe u {B_\nu^2} \vartheta_{n,\nu}=\sum_{i=0}^{\nu+\sigma(\nu)}
\vartheta_{\nu,i}\pe u {B_n^2}\delta_{i,n}=\vartheta_{\nu,n}
\pe u {B_n^2}.
$$
In particular, for $0\le \nu\le n-\sigma-1$, then $n\ge
\nu+\sigma+1\ge \nu+\sigma(\nu)+1$. Thus, we deduce
$\vartheta_{\nu,n}=0$. Hence $\vartheta_{n,\nu}=0$, for $0\le
\nu\le n-\sigma-1$.
\\
For $\nu=n-\sigma$, and $n\ge \sigma+t$, we obtain
$$
\ba{rl}
\pe u {B_{n-\sigma}^2}\vartheta_{n,n-\sigma}=&
\dst \pe u {\dq\big(\Phi(s)\dq B_{n-\sigma}(s-1)+
\Psi(s)B_{n-\sigma}(s-1)\big)}\\[0.3cm] & \dst
-\pe u {\dq\big(B_{n-\sigma}(s)
\dq \Phi(s-1)\big)B_n}=\sum_{nu=0}^{n+1}\tilde
\lambda_{n-\sigma,\nu} \pe u {B_n\dq B_\nu}\\[0.3cm]
=& \dst [n+1] \tilde \lambda_{n-\sigma,n+1} \pe u {B_n^2}.
\ea$$
But, from \refe{lambtildenp}, we get
$\vartheta_{n,n-\sigma}\ne 0$, either $n\ge \sigma+t+1$, or
$n=\sigma+t$ and $p\ge t-1$.
\\
As a consequence,
$$
\Phi(s)[\dq]^2 B_n(s-1)+\dq(\Psi(s)B_n(s-1))-B_n(s)[\dq]^2
\Phi(s-1)= \dst\sum_{\nu=n-\sigma}^{n+\sigma(n)}
{\vartheta}_{n,\nu}B_{\nu}(s),\ n\geq \sigma.
$$
{\bf (ii)$\Rightarrow$(i)}. From \refe{3.31}
$$
\ba{l}
\pe {\dq(\Phi(s-1)\dq u)+\big((\dq \Phi(s-1))-\Psi(s)\big)
\dq u}{B_n(s-1)}= 0,\quad n\ge \sigma+1,\\[0.3cm]
\pe {\dq(\Phi(s-1)\dq u)+\big((\dq \Phi(s-1))-\Psi(s)\big)
\dq u}{B_n(s-1)}= \pe u 1 \vartheta_{n,0}
,\quad n\le \sigma.\\[0.3cm]
\ea
$$
According to Lemma \ref{lem21}
$$\ba{rl}
\dq(\Phi(s-1)\dq u)+\big((\dq \Phi(s-1))-\Psi(s)\big) \dq u
=& \dst \sum_{n=0}^\sigma \frac{\pe u 1 \vartheta_{n,0}}
{\pe u {B_n^2}}B_n (\nabla u-u) \\[0.3cm]
= & \dst \sum_{n=0}^{\sigma(0)} \vartheta_{0,n}B_n(\nabla u-u).
\ea$$
Finally, a direct calculation yields
$$
\dq\big(\dq(\Phi u)-\Psi u\big)=0,
$$
then $\dq(\Phi u)-\Psi u=0$.
\\
Moreover, since $\sigma(n)=\sigma$ and $\vartheta_{n,n+\sigma}=
[n+\sigma+1]\tilde \lambda_{n,n+\sigma+1}\ne 0$, for $n\ge t+1$,
then $\tilde \lambda_{n,n+\sigma+1}\ne 0$, $n\ge t+1$.
The $q$-admissibility of the pair $(\Phi,\Psi)$ follows taking
into account the value of $\tilde \lambda_{n+\sigma(n)+1}$.
\edm \\
Our main result is the next one.
\bt \label{the3.4}
For any monic polynomial $\Phi$, $\deg \Phi=t$, and any
SMOP $\{B_n\}_{n\ge 0}$ with respect to $u$, the
following statements are equivalent.
\be
\item[\bf (i)] There exist a non-negative integer $\sigma$,
an integer $p\ge 1$, and an integer $r\ge \sigma+t+1$, with
$\sigma=\max(t-2,p-1)$, such that
\bq \label{aux7}
\sum_{\nu=n-\sigma}^{n+\sigma}\xi_{n,\nu}B_\nu(s)=
\sum_{\nu=n-t}^{n+\sigma}\varsigma_{n,\nu}B_\nu^{[1]}(s),
\eq
where $\, \xi_{n,n+\sigma}=\varsigma_{n,n+\sigma}=1$,
$n\ge \max(\sigma,t+1)$,
$\xi_{r,r-\sigma}\varsigma_{r,r-t}\ne 0$,
$$
\left\{\ba{ll}
\pe {\dq(\Phi u)} {B_m}=0,& p+1\le m\le 2
\sigma+t+1, \\[0.2cm] \pe {\dq(\Phi u)} {B_p}\ne 0,
\ea\right.
$$
and if $p=t-1$, then $\dst \lim_{q\uparrow 1}
\pe u {B_p^2}^{-1} \pe u {\Phi \dq B_p}\ne m$,
$m\in \XX N^*$ ($q$-admissibility condition).
\item[\bf (ii)] There exists a polynomial $\Psi$,
$\deg \Psi=p\ge 1$, such that
\bq \label{aux8}
\dq(\Phi u)=\Psi u,
\eq
where the pair $(\Phi,\Psi)$ is admissible.
\ee
\et
\bdm
{\bf (i)$\Rightarrow$ (ii)}. Let us consider the SMP
$\{\Xi_n\}_{n\ge 0}$ given by
$$
\ba{rl}
\Xi_{n+\sigma+1}(x)=& \dst \sum_{\nu=n-t}^{n+\sigma}
\frac{[n+\sigma+1]}{[\nu+1]}\varsigma_{n,\nu}B_{\nu+1}(x),
\quad n\ge \sigma+t+1, \\[0.3cm]
\Xi_n(x)= & B_n(x), \quad 0\le n\le 2\sigma+t+1.
\ea
$$
A direct calculation yields
$$
\dq \Xi_{n+\sigma+1}(s)=[n+\sigma+1]\sum_{\nu=n-\sigma}^{n
+\sigma} \xi_{n,\nu}B_\nu(s),\quad n\ge \sigma+t+1.
$$
Taking into account the linear functional $u$ is
quasi-definite, we get
$$
\pe {\dq(\Phi u)} {\Xi_{n+\sigma+1}}=-\pe u {\Phi \dq
\Xi_{n+\sigma+1}(s)}=-[n+\sigma+1]\!\!\sum_{\nu=n-
\sigma}^{n+\sigma} \!\!\xi_{n,\nu}\pe u {\Phi B_\nu}=0,
\quad n\ge \sigma+t+1.
$$
From the assumption and Lemma \ref{lem21}, if we denote
$\{w_n\}_{n\ge 0}$ the dual sequence of $\{\Xi_n\}_{n\ge 0}$,
then we get
\bq \label{aux9}
\dq (\Phi u)=\sum_{\nu=0}^{p} \pe {\dq(\Phi u)} {B_\nu} w_k.
\eq
Taking $\widehat t=\sigma+t$, $\widehat \sigma=\sigma+1$,
and $\widehat r=r+\sigma+1$, the polynomials $\Xi_n$ can
be rewritten as follows
$$
\ba{rl}
\Xi_n(x)= & \dst \sum_{\nu=n-\widehat t}^n\widetilde
\varsigma_{n,\nu} B_v(x),\quad n\ge \widehat \sigma+
\widehat t+1, \\[0.3cm]
\Xi_n(x)=& B_n(x),\quad 0\le n\le \widehat \sigma+\widehat t,
\ea
$$
where
$$
\ba{rl}
\widetilde \varsigma_{n,\nu}=& \dfrac{[n]}{[\nu]}\,
\varsigma_{n-\sigma-1, \nu-1},\quad n-\widehat t\le
\nu\le n, \ n\ge \sigma+\widehat t+1,\\[0.3cm]
\widetilde \varsigma_{\widehat r,\widehat r-\widehat t}=&
\dfrac{[r+\sigma+1]}{[r-t+1]}\varsigma_{r,r-t}\ne 0, \quad
\widehat r\ge 2\sigma+t+2\ge \widehat \sigma+\widehat t+1.
\ea
$$
From Lemma \ref{lem22}, $w_k=u_k=\pe u
{B_k^2}^{-1}B_k u$, $0\le k\le \widehat \sigma=\sigma+1$.
So, \refe{aux9} becomes
$$
\dq (\Phi u)=\sum_{\nu=1}^p \left(\frac{\pe {\dq(\Phi u)}
{B_\nu}} {\pe u {B_\nu^2}}B_\nu \right)u=\Psi u.
$$
Since $\pe {\dq(\Phi u)} {B_p}\ne 0$, then $\deg \Psi=p$.
\\
From the assumption, if $p=t-1$, then
$$
\lim_{q\uparrow 1}\frac 1 {[p]!}\left[\dq\right]^p\Psi(0)=
\lim_{q\uparrow 1} \frac{\pe {\dq(\Phi u)} {B_p}}
{\pe u {B_p^2}}=-\lim_{q\uparrow 1} \frac{\pe u
{\Phi\dq B_p}}{\pe u {B_p^2}}\ne -m, \quad m\in \XX N^*.
$$
Hence, the pair $(\Phi,\Psi)$ is admissible with associated
integer $\sigma$.
\\
{\bf (ii)$\Rightarrow$(i)}. From Lemma \ref{lem31}(iii),
there exists a polynomial $\Psi$, $\deg \Psi=p\ge 1$, such
that
\bq \label{aux10}
\Phi(s-1)\dq B_n(s-1)+\Psi(s)B_n(s-1)-B_n(s)\dq \Phi(s-1))=
\sum_{\nu=n-t+1}^{n+\sigma(n)+1}\tilde \lambda_{n,\nu}
B_\nu(s),\quad n\ge t,
\eq
where $\tilde \lambda_{n,n-t+1}\ne 0$, $n\ge t$,
$\sigma=\max(t-2,p-1)$, and the pair $(\Phi,\Psi)$ is
admissible.
\\
Taking $q$-differences in both hand sides of \refe{aux10}, we get
\bq \label{aux11}
\Phi(s)[\dq]^2 B_n(s-1)+\dq(\Psi(s)B_n(s-1))-B_n(s)[\dq]^2
\Phi(s-1)= \dst\sum_{\nu=n-t}^{n+\sigma(n)}
\zeta_{n,\nu}B_{\nu}^{[1]}(s),\ n\geq t,
\eq
where $\zeta_{n,\nu}=[\nu+1]\tilde \lambda_{n,\nu+1}$,
$0\le \nu\le n+\sigma(n)$, $n\ge t$.
\\
From \refe{3.31} and \refe{aux11}, we obtain \refe{aux7} where
$$
\ba{rl}
\xi_{n,\nu}=&\dfrac {\vartheta_{n,\nu}}{\vartheta_{n,n+\sigma}},
\quad n-\sigma\le \nu\le n+\sigma, \\[0.3cm]
\varsigma_{n,\nu}= & \dfrac{[\nu+1]\tilde \lambda_{n,\nu+1}}
{\vartheta_{n,n+\sigma}},\quad n-t\le \nu\le n+t,\\[0.3cm]
\xi_{n,n-\sigma}\varsigma_{n,n-t}= & \dfrac{[n-t+1]}
{\vartheta^2_{n,n+\sigma}}\,\vartheta_{n,n-\sigma}\tilde
\lambda_{n,n-t+1}\ne 0,\quad n\ge \sigma+t+1.
\ea
$$
Finally,
$$
\pe {\dq(\Phi u)} {B_n}=\pe u {\Psi B_n}=\left\{\ba{ll}
0, & p+1\le n\le 2\sigma+t+1, \\[0.2cm]
(\pe u {B_p^2}/[p]!)[\dq]^p\Psi(0)\ne 0, & n=p=\deg \Psi.
\ea\right.
$$
From the admissibility of the pair $(\Phi,\Psi)$, if $p=t-1$, then
$ \pe u {B_p^2}^{-1}\pe u {\Phi \dq B_p}\ne m, \ m\in \XX N^*$.
\edm
\section{The uniform lattice $x(s)=s$}
\noindent
As a direct consequence from the operator $L_{q,\omega}$ and the
$q$-linear lattice $x(s)$, we can recover the uniform lattice
setting $x(s)=(q^s-1)/(q-1)$ and taking limit $q\rightarrow 1$.
For instance, for $\Delta$-classical orthogonal polynomials the
structure relations \refe{fistre} and \refe{sestre} have been
studied in \cite{gamasa}.
\bt \label{theo31d} {\bf First Characterization of discrete
semiclassical polynomials}
\\
For a monic polynomial $\Phi$, $\deg \Phi=t$,
and any SMOP $\{B_n\}_{n\ge 0}$ with respect to $u$, the
following statements are equivalent.
\be
\item[(i)] There exist a non-negative integer $\sigma$, an
integer $p\ge 1$, and an integer $r\ge \sigma+t+1$, with
$\sigma=\max(t-2,p-1)$, such that
\bq \label{relteo1p}
\sum_{\nu=n-\sigma}^{n+t} \alpha_{n,\nu}B_\nu(s)=
\sum_{\nu=n-t}^{n+t} v_{n,\nu}B^{[1]}_\nu(s), \quad n\ge
\max(\sigma,t),
\eq
where $B_n^{[1]}(s):=(n+1)^{-1}\Delta B_{n+1}(s)$, $\alpha_{n,n+t}=
v_{n,n+t}=1$, $n\ge \max(\sigma,t)$, $\alpha_{r,r-\sigma}v_{r,r-t}\ne 0$,
$$
\langle\Delta(\Phi u),B_n\rangle=0, \ p+1\le n\le \sigma+
2t+1, \quad \langle \Delta(\Phi u),B_p\rangle\ne 0,
$$
and if $p=t-1$, then $\langle u,B_p^2\rangle^{-1}
\langle u,\Phi \Delta B_p\rangle\ne -m$, $m\in \XX N^*$.
\item[(ii)]There exists a polynomial $\Psi$,
$\deg \Psi=p\ge 1$, such that
$$
\Delta (\Phi u)=\Psi u,
$$
and the pair $(\Phi,\Psi)$ is admissible.
\ee
\et
\bt {\bf Second Characterization of discrete semiclassical polynomials}
For any monic polynomial $\Phi$, $\deg \Phi=t$, and any
SMOP $\{B_n\}_{n\ge 0}$ with respect to $u$, the
following statements are equivalent.
\be
\item[\bf (i)] There exist a non-negative integer $\sigma$,
an integer $p\ge 1$, and an integer $r\ge \sigma+t+1$, with
$\sigma=\max(t-2,p-1)$, such that
\bq \label{aux7p}
\sum_{\nu=n-\sigma}^{n+\sigma}\xi_{n,\nu}B_\nu(s)=
\sum_{\nu=n-t}^{n+\sigma}\varsigma_{n,\nu}B_\nu^{[1]}(s),
\eq
where $\xi_{n,n+\sigma}=\varsigma_{n,n+\sigma}=1$,
$n\ge \max(\sigma,t+1)$,
$\xi_{r,r-\sigma}\varsigma_{r,r-t}\ne 0$,
$$
\left\{\ba{ll}
\pe {\Delta(\Phi u)} {B_m}=0,& p+1\le m\le 2
\sigma+t+1, \\[0.2cm] \pe {\Delta(\Phi u)} {B_p}\ne 0,
\ea\right.
$$
and if $p=t-1$, then $\pe u {B_p^2}^{-1} \pe u {\Phi \Delta B_p}\ne m$,
$m\in \XX N^*$ (admissibility condition).
\item[\bf (ii)] There exists a polynomial $\Psi$,
$\deg \Psi=p\ge 1$, such that
\bq \label{aux8p}
\Delta(\Phi u)=\Psi u,
\eq
where the pair $(\Phi,\Psi)$ is admissible.
\ee
\et \noindent
The proofs are analogous to the original ones setting $\omega=1$,
and taking limit $q\uparrow 1$. Therefore $L_{q,1}\equiv \dq$ becomes
$\Delta$ and $[n]$ becomes $n$.
\bn $\Delta$-semiclassical linear functionals
have been studied in \cite{masa}.
\en
\section{examples}
\subsection{First example}\label{4.1}
Let $\{Q_n\}_{n\ge 0}$ be a SMOP that satisfies the
following relation
\bq \label{r2ex1}
(x(s+1)+v_{n,0})Q_n(s)=qQ_{n+1}^{[1]}(s)+\rho_n(s)Q_n^{[1]}(s),
\eq
where the lattice, $x(s)$, is $q$-linear, i.e.
$x(s+1)-qx(s)=\omega$,
$$
\ba{rl}
\rho_n= & \dfrac{q^{n+1}}{{\g C}}\dfrac{[n+1]}{\gamma_{n+1}},\
n\ge 1,\quad \rho_0=0,\\[0.4cm] v_{n,0}=&\dfrac{\gamma_{n+2}
\gamma_{n+1}}{q^{n}[n+2]}{\g C}+ \rho_n-q\beta_n-\omega,
\quad n\ge 0,
\ea
$$
and $\g C$ is a constant, being $\{\beta_n\}_{n\ge 0}$
and $\{\gamma_n\}_{n\ge 0}$ the coefficients of the TTRR
$$
xQ_n=Q_{n+1}+\beta_nQ_n+\gamma_nQ_{n-1}, \quad n\ge 1.
$$
Then, from the above TTRR and Theorem \ref{theo31},
we get $\{Q_n\}_{n\ge 0}$ is a sequence
of $q$-semiclassical orthogonal polynomials with respect to
the linear functional $v$, solution of the Pearson equation
\bq \label{eqv1}
\dq v=\Psi v,
\eq
of class $\sigma=1$, with $\Phi(x)=1$ and $\deg \Psi=2$.
\\
Then, it also satisfies the following relation
\bq \label{r1ex1}
Q_{n}^{[1]}(s)=Q_n(s)+\lambda_{n,n-1}Q_{n-1}(s),
\eq
where $\lambda_{n,n-1}=\dfrac{\gamma_{n+1}\gamma_n}
{q^n[n+1]}\g C$.\\
In fact, a straightforward calculation gives $\Psi(x)=-
\dfrac{\g C}q\, Q_2(x)-\dfrac 1 {\gamma_1}\,Q_1(x)$.
\bl 
Let $\{Q_n\}_{n\ge 0}$ be a SMOP with respect to the
linear functional $v$ satisfying \refe{eqv1}. Then
the sequence $\{Q_n\}_{n\ge 0}$ is not diagonal.
\el
\bdm
Assume $\{Q_n\}_{n\ge 0}$ is diagonal with respect to
$\phi$, with $\deg \phi=t$, and index $\sigma$. Then
from Corollary \ref{cor1}, $t/2\le \sigma\le t+2$
and we have the following diagonal relation
$$
\phi(s)Q_n(s)=\sum_{\nu=n-\sigma}^{n+t}\theta_{n,\nu}
Q_\nu^{[1]}(s),\quad \theta_{n,n-\sigma}\ne 0,\ n\ge \sigma.
$$
If we denote by $\{v_n\}_{n\ge 0}$ and $\{v_n^{[1]}\}_{n\ge 0}$
the dual sequences of $\{Q_n\}_{n\ge 0}$ and
$\{Q_n^{[1]}\}_{n\ge 0}$, respectively, then by Proposition \ref{lem16}
the last relation is equivalent to
\bq \label{aux4}
\phi v_n^{[1]}=k_n\Omega_{n+\sigma} v,\ n\ge 0,
\eq
where $k_n=\pe v {Q_{n+\sigma}^2}^{-1}\theta_{n+\sigma,n}$, and
$$
\Omega_{n+\sigma}(x)=\sum_{\nu=0}^{n+\sigma}\frac{\theta_{\nu,n}}
{\theta_{n+\sigma,n}}\frac{\pe v {Q_{n+\sigma}^2}}{\pe v
{Q_{\nu}^2}} Q_\nu(x),\quad n\ge 0.
$$
It is clear that $v$ satisfies an infinite number of
relations as \refe{aux4}. Indeed, by multiplying both
hand sides of \refe{aux4} by a monic polynomial, we get
another diagonal relation.
\\
For this reason, we will assume $t=\deg \phi$ is the
minimum non-negative integer such that $v$ satisfies diagonal
relations as \refe{aux4}, i.e. the Eq. \refe{aux4} cannot
be simplified.
\\
Notice that $t\ge 1$. Indeed, if we suppose that $t=0$,
then $0\le \sigma\le 2$ and we recover the first structure
relation characterizing $q$-classical sequences.
This contradicts the fact that the sequence
$\{Q_n\}_{n\ge 0}$ is $q$-semiclassical of class one.
\\
Consequently, since $t\ge 1$ then $\sigma\ge 1$.
Taking $q$-differences in both hand sides of \refe{aux4}
and using \refe{2.1}, from \refe{eqv1} and
$\dq v_n^{[1]}=-[n+1]v_{n+1}$, we obtain
\bq \label{aux5}
\widetilde \phi v_n^{[1]}=k_n\psi_n v, \quad n\ge 0,
\eq
where
$$ \ba{rl}
\widetilde \phi(s)=& [t]^{-1}\dq \phi(s),\\[0.3cm]
\psi_n(s)=& [t]^{-1}\big(\Omega_{n+\sigma}(s+1)
\Psi(s)+\dq \Omega_{n+\sigma}(s)
+d_n\phi(s+1)Q_{n+1}(s)\big), \ n\ge 0,\\[0.3cm]
d_n=& [n+1](\pe v {Q_{n+1}^2} k_n)^{-1}, \ n\ge 0.
\ea$$
Notice that the polynomial $\widetilde \phi$ is monic
with $\deg \widetilde \phi=t-1$.
\\
Moreover, taking into account $u$ is a quasi-definite
linear functional, combining \refe{aux4} and \refe{aux5}
we obtain $\widetilde \phi(x)\Omega_{n+\sigma}(x)=
\phi(x)\psi_n(x)$, and analyzing
the highest degree of this relation, we get $\psi_n$
is a monic polynomial with $\deg \psi_{n}=n+\sigma-1$.
But, this contradicts the fact that $t=\deg \phi$ is
the minimum nonnegative integer such that $v$ satisfies
diagonal relations as \refe{aux4}.
\edm
\subsection{The $q$-Freud type polynomials}
Let $\{P_n\}_{n\ge 0}$ be a SMOP with respect to a linear
functional $u$ such that $(u)_0=\pe u 1=1$ and the following relation
\bq \label{defex1}
\dq P_n(s)=[n]P_{n-1}(s)+a_n P_{n-3}(s),\quad n\ge 2,
\eq
holds, where $P_{-1}\equiv 0$, $P_0\equiv 1$, and $P_1(x)=x$,
being $x\equiv x(s)=q^s$, i.e. $\omega=0$.
\\
We know that this family satisfies a TTRR, i.e. there
exist two sequences of complex numbers $\{b_n\}_n$ and
$\{c_n\}_n$, $c_n\ne 0$, such that
$$
x P_n=P_{n+1}+b_n P_n+c_nP_{n-1}, \quad n\ge 1.
$$
Furthermore, from a direct calculation we get
$a_n=K(q) q^{-n}c_{n}c_{n-1}c_{n-2}$, $n\ge 2$. In fact,
the parameters $c_n$ satisfy the non-linear recurrence
relation
$$
q[n]c_{n-1}+K(q)q^{-n+1}c_nc_{n-1}c_{n-2}=[n-1]c_n+
K(q)q^{-n-1}c_{n+1}c_nc_{n-1}, \quad n\ge 1,
$$
with $c_0=0$, $c_1=-P_2(0)\ne 0$, and
$\dst \lim_{q\uparrow 1}K(q)=4$.
\\
Moreover, from Proposition \ref{lem31} we deduce that
$\Phi\equiv 1$ and thus $\sigma=2$. As a consequence $\Psi$ is a
polynomial of degree 3. In other words, $u$ is a
$q$-semiclassical linear functional of class 2, i.e. $u$
satisfies the following distributional equation
\bq \label{qfreud}
\dq u=\Psi u, \quad \deg \Psi=3.
\eq
\bl \label{lem1}
$\Psi(x)=-K(q)q^{-3}P_3(x)-c_1^{-1}P_1(x)$.
\el
\noindent
So, \refe{qfreud} is the $q$-analog of the
Pearson equation for the Freud case.
\\
\bdm
From our hypothesis $\Psi$ is a polynomial of degree 3, so
$\Psi(x)=e_0P_0+e_1P_1+e_2P_2+e_3P_3$.
Then, taking into account $d_n^2=c_nd_{n-1}^2$,
$n\ge 1$, and the value of $a_n$, $n\ge 3$, we get
$$
\ba{l}
e_0d_0^2=e_0\pe u {P_0^2} =\pe {\Psi u} {P_0}=
-\pe u {\dq P_0}=0, \\
e_1d_1^2=e_1\pe u {P_1^2} =\pe {\Psi u} {P_1}=
-\pe u {\dq P_1}=-1, \\
e_2d_2^2=e_2\pe u {P_2^2} =\pe {\Psi u} {P_2}=
-\pe u {\dq P_2}
\stackrel{\refe{defex1}}=-\pe u {[2]P_1}=0, \\
e_3d_3^2=e_3\pe u {P_3^2} =\pe {\Psi u} {P_3}=
-\pe u {\dq P_3}
\stackrel{\refe{defex1}}=-\pe u {[3]P_2+a_3P_0}=-a_3.
\ea
$$
\edm
\\
From Theorem \ref{the3.4}, we can
write the second structure relation as follows
\bq \label{qfressr}
B_{n+2}+\xi_{n,n+1}B_{n+1}+\xi_{n,n}B_n+\xi_{n,n-1}B_{n-1}
+\xi_{n,n-2}B_{n-2}=B_{n+2}^{[1]}+\varsigma_{n,n+1} B_{n+1}^{[1]}+
\varsigma_{n,n}B_{n}^{[1]}.
\eq
Using \refe{defex1} we get
$$ \ba{ll}
\xi_{n,n+1}=\varsigma_{n,n+1},& \quad \xi_{n,n}=[n+3]^{-1}a_{n+3}+
\varsigma_{n,n},\\ \xi_{n,n-1}=[n+2]^{-1}\varsigma_{n,n+1}a_{n+2},&
\quad \xi_{n,n-2}=[n+1]^{-1}\varsigma_{n,n}a_{n+1}. \ea
$$
Moreover, combining both structure relations if $P_n(x)=\sum_{j=0}^n
\lambda_{n,j}\,x^{n-j}$, then $\lambda_{n,2k+1}=0$ for
nonnegative integers $n, \, k$ such that $0\le k\le (n-1)/2$, and
$$
\lambda_{n,0}=1,\quad \lambda_{n,2k+2}=\frac{[n]c_{n-1}\lambda_{n-2,
2k}+a_n\lambda_{n-3,2k}}{[n-2k-2]-[n]},\quad 1\le k\le n/2.
$$
In fact, with these values, we obtain $c_n=\lambda_{n,2}-
\lambda_{n+1,2}$, $b_n=\lambda_{n,1}-\lambda_{n+1}=0$, and
$\xi_{n,n+1}=\xi_{n,n-1}=\varsigma_{n,n+1}=0$,
$n\ge 0$. Hence, we can rewrite \refe{qfressr} as
\bq \label{qfressrp}
(x^2+\widetilde{v}_{n,0})B_n=B_{n+2}^{[1]}+\widetilde{\rho}_nB_n^{[1]},
\eq
where $\widetilde{v}_{n,0}=\dfrac{a_{n+3}}{[n+3]}+\dfrac{q^{n+1}[n+1]}
{K(q)c_{n+1}}-c_{n+1}-c_n$, and $\widetilde{\rho}_n=
\dfrac{q^{n+1}[n+1]}{K(q)c_{n+1}}$.
\\
\bl
The moments of the linear functional $u$, $\{(u)_n\}_{n\ge 0}$,
satisfy the following relation
\bq \label{momrel}
[n+1](u)_{n}=K(q)q^{-3}(u)_{n+4}+\left(\frac 1 {c_1}-\dfrac{[3]c_2+a_3}{q(1+q)}
\right)(u)_{n+2},\quad n\ge 0,
\eq
where $(u)_0=1$.
\el \noindent
Therefore, taking into account that $(u)_1=(u)_3=0$, we can deduce $u$ is a
symmetric linear functional, i.e. $(u)_{2n+1}=\pe u {x^{2n+1}}=0$,
$n\ge 0$.
\\[0.4cm]
{\bf Acknowledgements.} This work has been supported by
Direcci\'{o}n General de Investigaci\'{o}n (Ministerio de Educaci\'{o}n
y Ciencia) of Spain, grant BFM03-06335-C03-02 and INTAS
research network NeCCA 03-51-6637.

\end{document}